# Gibbs fragmentation trees

PETER MCCULLAGH[1], JIM PITMAN[2] and MATTHIAS WINKEL[3]

[1]*Department of Statistics, University of Chicago, 5734 University Ave, Chicago, IL 60637, USA.
E-mail: [pmcc@galton.uchicago.edu](pmcc@galton.uchicago.edu)*

[2]*Statistics Department, 367 Evans Hall # 3860, University of California, Berkeley, CA 94720-3860, USA. E-mail: [pitman@stat.berkeley.edu](pitman@stat.berkeley.edu)*

[3]*Department of Statistics, University of Oxford, 1 South Parks Road, Oxford OX1 3TG, UK.
E-mail: [winkel@stats.ox.ac.uk](winkel@stats.ox.ac.uk)*

We study fragmentation trees of Gibbs type. In the binary case, we identify the most general Gibbs-type fragmentation tree with Aldous' beta-splitting model, which has an extended parameter range $\beta > -2$ with respect to the beta$(\beta+1, \beta+1)$ probability distributions on which it is based. In the multifurcating case, we show that Gibbs fragmentation trees are associated with the two-parameter Poisson–Dirichlet models for exchangeable random partitions of $\mathbb{N}$, with an extended parameter range $0 \leq \alpha \leq 1$, $\theta \geq -2\alpha$ and $\alpha < 0$, $\theta = -m\alpha$, $m \in \mathbb{N}$.

## 1. Introduction

We are interested in various models for random trees associated with processes of recursive partitioning of a finite or infinite set, known as *fragmentation processes* [2, 4, 9]. We start by introducing a convenient formalism for the kind of combinatorial trees arising naturally in this context [16, 18]. Let $\#B$ be the number of elements in the finite non-empty set $B$. Following standard terminology, a *partition of* $B$ is a collection

$$\pi_B = \{B_1, \ldots, B_k\}$$

of non-empty disjoint subsets of $B$ whose union is $B$. To introduce a new terminology convenient for our purpose, we make the following recursive definition. A *fragmentation of* $B$ (sometimes called a *hierarchy* or a *total partition*) is a collection $\mathbf{t}_B$ of non-empty subsets of $B$ such that

(i) $B \in \mathbf{t}_B$;

(ii) if $\#B \geq 2$ then, there is a partition $\pi_B$ of $B$ into $k$ parts, $B_1, \ldots, B_k$, called the *children of $B$*, for some $k \geq 2$, with

$$\mathbf{t}_B = \{B\} \cup \mathbf{t}_{B_1} \cup \cdots \cup \mathbf{t}_{B_k}, \tag{1}$$

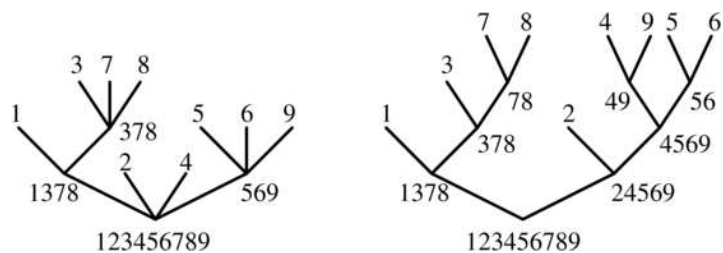

**Figure 1.** Two fragmentations of [9] graphically represented as trees labeled by subsets of [9].

where $\mathbf{t}_{B_i}$ is a fragmentation of $B_i$ for each $1 \leq i \leq k$.

Necessarily, $B_i \in \mathbf{t}_B$, each child $B_i$ of $B$ with $\#B_i \geq 2$ has further children, and so on, until the set $B$ is broken down into singletons. We use the same notation $\mathbf{t}_B$ both

- for such a collection of subsets of $B$, and
- for the tree whose vertices are these subsets of $B$ and whose edges are defined by the parent/child relation determined by the fragmentation.

To emphasize the tree structure, we may call $\mathbf{t}_B$ a *fragmentation tree*. Thus, $B$ is the root of $\mathbf{t}_B$ and each singleton subset of $B$ is a leaf of $\mathbf{t}_B$ (see Figure 1 – here $[9] = \{1, \ldots, 9\}$; we also put $[n] = \{1, \ldots, n\}$). We denote by $\mathbb{T}_B$ the collection of all fragmentations of $B$. A fragmentation $\mathbf{t}_B \in \mathbb{T}_B$ is called *binary* if every $A \in \mathbf{t}_B$ has either 0 or 2 children. We denote by $\mathbb{B}_B \subseteq \mathbb{T}_B$ the collection of binary fragmentations of $B$.

For each non-empty subset $A$ of $B$, the *restriction to $A$ of $\mathbf{t}_B$*, denoted $\mathbf{t}_{A,B}$, is the fragmentation tree whose root is $A$, whose leaves are the singleton subsets of $A$ and whose tree structure is defined by restriction of $\mathbf{t}_B$. That is, $\mathbf{t}_{A,B}$ is the fragmentation $\{C \cap A : C \cap A \neq \varnothing, C \in \mathbf{t}_B\} \in \mathbb{T}_A$, corresponding to a *reduced subtree*, as discussed by Aldous [1].

Given a rooted combinatorial tree with no single-child vertices and whose leaves are labeled by a finite set $B$, there is a corresponding fragmentation $\mathbf{t}_B$, where each vertex of the combinatorial tree is associated with the set of leaves in the subtree above that vertex. So the fragmentations defined here provide a convenient way to label the vertices of a combinatorial tree and to encode the tree structure in the labeling.

A *random fragmentation model* is an assignment, for each finite subset $B$ of $\mathbb{N}$, of a probability distribution on $\mathbb{T}_B$ for a random fragmentation $T_B$ of $B$. We assume throughout this paper that the model is *exchangeable*, meaning that the distribution of $T_B$ is invariant under the obvious action of permutations of $B$ on fragmentations of $B$. The distribution of $\Pi_B$, the partition of $B$ generated by the branching of $T_B$ at its root, is then of the form

$$\mathbb{P}(\Pi_B = \{B_1, \ldots, B_k\}) = p(\#B_1, \ldots, \#B_k) \qquad (2)$$

for all partitions $\{B_1, \ldots, B_k\}$ with $k \geq 2$ blocks and some symmetric function $p$ of compositions of positive integers, called a *splitting probability rule*. The model is called

- *consistent* if for every $A \subset B$, the restricted tree $T_{A,B}$ is distributed like $T_A$;
- *Markovian* if, given $\Pi_B = \{B_1, \ldots, B_k\}$, the $k$ restricted trees $T_{B_1,B}, \ldots, T_{B_k,B}$ are independent and distributed as $T_{B_1}, \ldots, T_{B_k}$;
- *binary* if $T_B$ is a binary tree with probability one, for every $B$.

Aldous [2] initiated the study of consistent Markovian binary trees as models for neutral evolutionary trees. He observed parallels between these models and Kingman's theory of exchangeable random partitions of $\mathbb{N}$, and posed the problem of characterizing these models analogously to known characterizations of the Ewens sampling formula for random partitions. In [9], we showed how consistent Markovian trees arise naturally in Bertoin's theory of homogeneous fragmentation processes [4] and deduced from Bertoin's theory a general integral representation for the splitting rule of a Markovian fragmentation model.

To briefly review these developments in the binary case, the distribution of a Markovian binary fragmentation $T_B$ is determined by a splitting rule $p$, which is a symmetric function $p$ of pairs of positive integers $(i,j)$, according to the following formula for the probability of a given tree $\mathbf{t} \in \mathbb{B}_B$:

$$\mathbb{P}(T_B = \mathbf{t}) = \prod_{A \in \mathbf{t}: \#A \geq 2} p(\#A_1, \#A_2), \tag{3}$$

where $A_1$ and $A_2$ denote the two children of $A$ in the tree $T_B$.

The following proposition collects some known results.

**Proposition 1.** (i) *Every non-negative symmetric function $p$ subject to normalization conditions*

$$\sum_{k=1}^{n-1} \binom{n-1}{k-1} p(k, n-k) = 1 \qquad \text{for all } n \geq 2$$

*defines a Markovian binary fragmentation model.*

(ii) *A splitting rule $p$ gives rise to a consistent Markovian binary fragmentation if and only if*

$$p(i,j) = p(i+1,j) + p(i,j+1) + p(i+j,1)p(i,j) \qquad \text{for all } i,j \geq 1. \tag{4}$$

(iii) *Every consistent splitting rule admits an integral representation*

$$p(i,j) = \frac{1}{Z(i+j)} \left( \int_{(0,1)} x^i (1-x)^j \nu(\mathrm{d}x) + c \mathbf{1}_{\{i=1 \text{ or } j=1\}} \right) \qquad \text{for all } i,j \geq 1, \tag{5}$$

*with characteristics $c \geq 0$ and $\nu$ a symmetric measure on $(0,1)$ with $\int_{(0,1)} x(1-x)\nu(\mathrm{d}x) < \infty$, and $Z(n)$ a sequence of normalization constants.*

**Proof.** (i) is elementary. For (ii), Ford [6], Proposition 41, gave a characterizaton of consistency for models of unlabeled trees which is easily shown to be equivalent to the condition stated here. The interpretation (and sketch of proof) of this condition is that for $B = C \cup \{k\}$ (with $k \notin C$), the vertex $C$ of $T_C$ splits into a particular partition of sizes $i$ and $j$ if and only if $T_B$ splits into that partition with $k$ added to one or the other block, or if $T_B$ first splits into $C$ and $\{k\}$ and then $C$ splits further into that partition of sizes $i$ and $j$. (iii) is directly read from [9]. □

Aldous [2] studied in some detail the *beta-splitting model* which arises as the particular case of (5) with characteristics $c = 0$ and

$$\nu(\mathrm{d}x) = x^\beta (1-x)^\beta \mathrm{d}x \quad \text{for } \beta \in (-2, \infty) \quad \text{and} \quad \nu(\mathrm{d}x) = \delta_{1/2}(\mathrm{d}x) \quad \text{for } \beta = \infty. \quad (6)$$

Aldous posed the problem of characterizing this model among all consistent binary Markov models. The main focus of this paper is the following result.

**Theorem 2.** *Aldous' beta-splitting models for $\beta \in (-2, \infty]$ are the only consistent Markovian binary fragmentations with splitting rule of the form*

$$p(i,j) = \frac{w(i)w(j)}{Z(i+j)} \quad \text{for all } i, j \geq 1, \quad (7)$$

*for some sequence of weights $w(j) \geq 0$, $j \geq 1$, and normalization constants $Z(n)$, $n \geq 2$.*

As a corollary, we extract a statement purely about measures on $(0,1)$.

**Corollary 3.** *Every symmetric measure $\nu$ on $(0,1)$ with $\int_{(0,1)} x(1-x)\nu(\mathrm{d}x) < \infty$, whose moments factorize into the form*

$$\int_{(0,1)} x^i (1-x)^j \nu(\mathrm{d}x) = w(i)w(j) \quad \text{for all } i, j \geq 1$$

*for some $w(i) \geq 0$, $i \geq 1$, is a multiple of one of Aldous' beta-splitting measures (6).*

In particular, this characterizes the symmetric beta distributions among probability measures on $(0,1)$.

Berestycki and Pitman [3] encountered a different one-dimensional class of Gibbs splitting rules in the study of fragmentation processes related to the affine coalescent. These are not consistent, but the Gibbs fragmentations are naturally embedded in continuous time.

The rest of this paper is organized as follows. Section 2 offers an alternative characterization of what we call *binary Gibbs models*, meaning models with splitting rule of the form (7), without assuming consistency. Theorem 2 is then proved in Section 3. In Section 4, we discuss growth procedures and embedding in continuous time for the consistent case. Section 5 gives a generalization of the Gibbs results to multifurcating trees.

## 2. Characterization of binary Gibbs fragmentations

The Gibbs model (7) is overparameterized: if we multiply $w(k)$, $k \geq 1$, by $ab^k$ (and then $Z(m)$, $m \geq 2$, by $a^2 b^m$), the model remains unchanged. Note, further, that neither $w(1) = 0$ nor $w(2) = 0$ is possible since then (7) does not define a probability function for $n = i + j = 3$. Hence, we may assume $w(1) = 1$ and $w(2) = 1$. It is now easy to see that for any two different such sequences, the models are different. Note that the following result does not assume a consistent model.

**Proposition 4.** *The following two conditions on a collection of random binary fragmentations $T_B$ indexed by finite subsets $B$ of $\mathbb{N}$ are equivalent:*

(i) *$T_B$ is for each $B$ an exchangeable Markovian binary fragmentation with splitting rule of the Gibbs form (7) for some sequence of weights $w(j) > 0$, $j \geq 1$, and normalization constants $Z(n)$, $n \geq 2$;*

(ii) *for each $B$, the probability distribution of $T_B$ is of the form*

$$\mathbb{P}(T_B = \mathbf{t}) = \frac{1}{w(\#B)} \prod_{A \in \mathbf{t}} \psi(\#A) \qquad \text{for all } \mathbf{t} \in \mathbb{B}_B, \tag{8}$$

*for some sequence of weights $\psi(j) > 0$, $j \geq 1$, and normalisation constants $w(n)$, $n \geq 1$.*

*More precisely, if* (i) *holds with $w(1) = 1$, then* (ii) *holds for the same sequence $w$ with*

$$\psi(1) = 1 \quad \text{and} \quad \psi(k) = w(k)/Z(k), \qquad k \geq 2. \tag{9}$$

*Conversely, if* (ii) *holds for some sequence $\psi$ with $\psi(1) = 1$, then* (i) *holds for the sequence $w(n)$, $n \geq 1$, determined by (8); in particular, $w(1) = 1$.*

**Proof.** Given a Gibbs model with $w(1) = 1$, we can combine (3) and (7) to get, for all $\mathbf{t} \in \mathbb{B}_B$,

$$\mathbb{P}(T_B = \mathbf{t}) = \prod_{A \in \mathbf{t}: \#A \geq 2} \frac{w(\#A_1) w(\#A_2)}{Z(\#A)} = \frac{1}{w(\#B)} \prod_{A \in \mathbf{t}: \#A \geq 2} \frac{w(\#A)}{Z(\#A)}.$$

If we make the substitution (9), we can read off $w(n)$ as the correct normalization constant and (8) follows, with $\psi(1) = 1$.

On the other hand, (8) determines the sequence $w(n)$, $n \geq 1$, as

$$w(n) = \sum_{\mathbf{t} \in \mathbb{B}_{[n]}} \prod_{A \in \mathbf{t}} \psi(\#A).$$

Note, in particular, that $w(1) = \psi(1)$. We can express the normalization constants in the Gibbs model (7) by the formula

$$Z(m) = \sum_{k=1}^{m-1} \binom{m-1}{k-1} w(k) w(m-k) \tag{10}$$

$$= \sum_{k=1}^{m-1} \binom{m-1}{k-1} \left( \sum_{\mathbf{t}_1 \in \mathbb{B}_{[k]}} \prod_{A \in \mathbf{t}_1} \psi(\#A) \right) \left( \sum_{\mathbf{t}_2 \in \mathbb{B}_{[m-k]}} \prod_{A \in \mathbf{t}_2} \psi(\#A) \right)$$

$$= \sum_{\mathbf{t} \in \mathbb{B}_{[m]}} \prod_{A \in \mathbf{t}: A \neq [m]} \psi(\#A) = w(m)/\psi(m),$$

as in (9). By application of the previous implication from (i) to (ii), formula (8) gives the distribution of the Gibbs model derived from this weight sequence $w(n)$ and the conclusion follows. □

Note that the normalization constant $Z(m)$ in the Gibbs splitting rule (7) model and given in (10) is a partial Bell polynomial in $w(1), w(2), \ldots$ (see [15] for more applications of Bell polynomials), whereas the normalization constant $w(n)$ in the Gibbs tree formula (8) is a polynomial in $\psi(1), \psi(2), \ldots$ of a much a more complicated form. The normalization constant in (8) is

$$w(n) = \sum_{\mathbf{t} \in \mathbb{B}_{[n]}} \prod_{A \in \mathbf{t}} \psi(\#A).$$

In an attempt to study this polynomial in $\psi(1), \psi(2), \ldots$, we introduce the *signature* $\sigma_\mathbf{t} : [n] \to \mathbb{N}$ of a tree $\mathbf{t} \in \mathbb{B}_{[n]}$ by

$$\sigma_\mathbf{t}(j) = \#\{A \in \mathbf{t} : \#A = j\}, \qquad j = 1, \ldots, n.$$

Note that $\mathbb{P}(T_n = \mathbf{t})$ depends on $\mathbf{t}$ only via $\sigma_\mathbf{t}$, that is, $\sigma_\mathbf{t}$ is a sufficient statistic for the Gibbs probabilities (8). Denote the set of signatures by $\mathrm{Sig}_n = \{\sigma_\mathbf{t} : \mathbf{t} \in \mathbb{B}_{[n]}\}$. The inductive definition of $\mathbb{B}_{[n]}$ yields

$$\mathrm{Sig}_n = \{\sigma^{(1)} + \sigma^{(2)} + 1_n : \sigma^{(1)} \in \mathrm{Sig}_{n_1}, \sigma^{(2)} \in \mathrm{Sig}_{n_2}, n_1 + n_2 = n\},$$

where $1_n(j) = 1$ if $j = n$, $1_n(j) = 0$ otherwise. The coefficients $Q_\sigma$ in $w(n)$, when expanded as a polynomial in $\psi(1), \psi(2), \ldots$, are numbers of fragmentations with the same signature $\sigma \in \mathrm{Sig}_n$:

$$w(n) = \sum_{\sigma \in \mathrm{Sig}_n} Q_\sigma \psi^\sigma, \qquad \text{where } \psi^\sigma = \prod_{j=1}^n \psi(j)^{\sigma(j)}.$$

Let us associate with each fragmentation $\mathbf{t} \in \mathbb{B}_{[n]}$ its tree shape (combinatorial tree without labels) $\mathbf{t}^\circ$ and denote by $\mathbb{B}_n^\circ$ the collection of shapes of binary trees with $n$ leaves. Clearly, two fragmentations with the same tree shape have the same signature, so we can define $\sigma(\mathbf{t}^\circ)$ in the obvious way. For $n \leq 8$ (and many larger trees), direct enumeration shows that the tree shape $\mathbf{t}^\circ \in \mathbb{B}_n^\circ$ is uniquely determined by its signature $\sigma$, and $Q_\sigma$ is just the number $q(\mathbf{t}^\circ)$ of different labelings. For $n \geq 9$, this is false: there are two tree shapes with signature $(9, 3, 1, 2, 1, 0, 0, 0, 1)$; see Figure 2. If we denote by $\mathcal{I}_\sigma^\circ \subseteq \mathbb{B}_n^\circ$ the set of tree shapes with signature $\sigma$, then $Q_\sigma = \sum_{\mathbf{t}^\circ \in \mathcal{I}_\sigma^\circ} q(\mathbf{t}^\circ)$. The remaining combinatorial problem is therefore to study $\mathcal{I}_\sigma^\circ$ and $q(\mathbf{t}^\circ)$. We have not been able to solve this problem. The preprint version [12] of the present paper includes an Appendix with a partial study: see also Corollary 2.4.3 of [17].

## 3. Consistent binary Gibbs rules

The statement of Theorem 2 specifies Aldous' [2] beta-splitting models by their integral representation (5). Observe that the moment formula for beta distributions easily gives

$$p(i,j) = \frac{1}{Z(i+j)} \int_0^1 x^{i+\beta}(1-x)^{j+\beta}\,dx$$
$$= \frac{\Gamma(i+\beta+1)\Gamma(j+\beta+1)}{R(i+j)} \qquad \text{for all } i,j \geq 1, \tag{11}$$

for normalization constants $R(n) = Z(n)\Gamma(n+2\beta+2)$, $n \geq 2$. This is for $\beta \in (-2, \infty)$. For $\beta = \infty$, we simply get $p(i,j) = 1/R(i+j)$ for all $i,j \geq 1$, where $R(n) = Z(n)2^n$, $n \geq 2$.

**Proof of Theorem 2.** We start from a general Gibbs model (7) with $w(1) = 1$ and follow [7], Section 2 closely, where a similar characterization is derived in a partition rather than a tree context. Let the Gibbs model be consistent. This immediately implies that $w(j) > 0$ for all $j \geq 1$. The consistency criterion (4) in terms of $W_j = w(j+1)/w(j)$ now gives

$$W_i + W_j = \frac{Z(i+j+1) - w(i+j)}{Z(i+j)} \qquad \text{for all } i,j \geq 1. \tag{12}$$

The right-hand side is a function of $i+j$, so $W_{j+1} - W_j$ is constant and hence $W_j = a + bj$ for some $b \geq 0$ and $a > -b$. Now, either $b = 0$ (excluded for the time being) or

$$w(j) = W_1 \cdots W_{j-1} = \prod_{q=1}^{j-1}(a+bq)$$

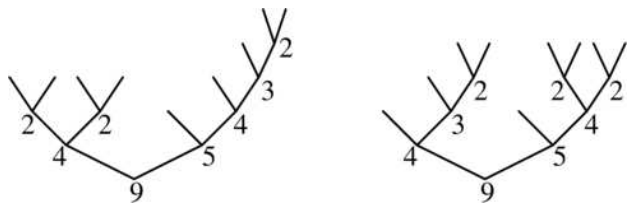

**Figure 2.** Two tree shapes with the same signature (here marked by subtree sizes).

$$= b^{j-1} \prod_{q=1}^{j-1} \left(\frac{a}{b} + q\right) = b^{j-1} \frac{\Gamma(a/b+j)}{\Gamma(a/b+1)}.$$

and, hence, reparameterizing by $\beta = a/b - 1 \in (-2, \infty)$ and pushing $b^{i+j-2}$ into the normalization constant $d_{i+j} = b^{i+j-2}/Z(i+j)$, we have

$$p(i,j) = \frac{w(i)w(j)}{Z(i+j)} = d_{i+j} \frac{\Gamma(i+1+\beta)}{\Gamma(2+\beta)} \frac{\Gamma(j+1+\beta)}{\Gamma(2+\beta)}.$$

The case $b = 0$ is the limiting case $\beta = \infty$, when, clearly, $w(j) \equiv 1$ (now pushing $a^{i+j-2}$ into the normalization constant).

These are precisely Aldous' beta-splitting models, as in (11).  □

While we identified the boundary case $\beta = \infty$ as being of Gibbs type, the boundary case $\beta = -2$ is not of Gibbs type, although it can still be made precise as a Markovian fragmentation model with characteristics $c > 0$ and $\nu = 0$ (pure erosion): $p(i,j) = 0$ unless $i = 1$ or $j = 1$, so the Markovian fragmentations $T_n$ are combs, where all $n - 1$ branching vertices are lined up in a single spine.

In the proof of the theorem, we obtained as parameterization for the Gibbs models (7),

$$w(j) = \frac{\Gamma(j+1+\beta)}{\Gamma(2+\beta)}, \qquad j \geq 1, \tag{13}$$

for some $\beta \in (-2, \infty)$, or $w(j) \equiv 1$ for $\beta = \infty$. Note that the simple convention $w(2) = 1$ from Section 2 is not useful here. We can now still deduce the parameterization (8) by Proposition 4, in principle. However, since $\psi(k) = w(k)/Z(k)$ involves partial Bell polynomials $Z(k)$ in $w(1), w(2), \ldots$, this is less explicit in terms of $\beta$ than the parameterization (7).

$$\psi(2) = 2 + \beta, \qquad \psi(3) = \frac{3+\beta}{3}, \qquad \psi(4) = \frac{(3+\beta)(4+\beta)}{18+7\beta}, \ldots.$$

Special cases that have been studied in various biology and computer science contexts (see Aldous [2] for a review) include the following: $\beta = -3/2, -1, 0, \infty$. In these cases, we can explicitly calculate the Gibbs parameters in (7) and (8) and the normalisation constants.

If $\beta = -3/2$, we can take $\psi(n) \equiv 1$ and $T_B$ is *uniformly distributed*: if $\#B = n$, then $\mathbb{P}(T_B = \mathbf{t}) = 2^{n-1}(n-1)!/(2n-2)!$, $\mathbf{t} \in \mathbb{B}_B$. The asymptotics of uniform trees lead to Aldous' Brownian CRT [1]; see also [15], Section 6.3. Table 1 uses a different parameterization via the convenient relations (9) and (13).

The case $\beta = -1$ is the limiting conditional distribution in the Ewens family as the Ewens parameter $\lambda \to 0$, conditional on the occurrence of a split. The $\beta = 0$ case is known as the *Yule model* and $\beta = \infty$ as the *symmetric binary trie* (see Aldous [2]). Continuum tree limits of the beta-splitting model for $\beta \in (-2, -1)$ are described in [9].

The normalization that leads to a compact limit tree is here $T_{[n]}/n^{-\beta-1}$, where $T_{[n]}$ is represented as a metric tree with unit edge lengths and the scaling $T_{[n]}/n^{-\beta-1}$ refers to scaling of edge lengths. Aldous [2] studies weaker asymptotic properties for average distance from a leaf to the root, also for $\beta \geq -1$, where growth is logarithmic.

## 4. Growth rules and embedding in continuous time

In [9], we study the consistently growing sequence $T_n$, $n \geq 1$, where $T_n := T_{[n]} = T_{[n],[n+1]}$ is the restriction of $T_{n+1}$ to $[n]$ for all $n \geq 1$, in a general context of consistent Markovian multifurcating fragmentation models. The integral representation (5) stems from an association with Bertoin's theory of homogeneous fragmentation processes in continuous time [4]. Let us here look at the binary case in general and Gibbs fragmentations in particular.

Consider the distribution of $T_{n+1}$, given $T_n$. The tree $T_{n+1}$ has a vertex $A \cup \{n+1\}$ with children $\{n+1\}$ and $A \in T_n$. We say that $n+1$ has been *attached below* $A$. In passing from $T_n$ to $T_{n+1}$, leaf $n+1$ can be attached below any vertex $A$ of $T_n$ (including $[n]$ and all leaf nodes). Note that to construct $T_{n+1}$ from $T_n$, $n+1$ is also added as an element to all vertices on the path from $[n]$ to $A$. Vertex $A \in T_n$ is special in that both $A$ and $A \cup \{n+1\}$ are in $T_{n+1}$.

Fix a vertex $A$ of $\mathbf{t} \in \mathbb{B}_{[n]}$ and consider the conditional probability, given $T_n = \mathbf{t}$, of $n+1$ being attached below $A$. This is the ratio of two probabilities of the form (3) in which many common factors cancel so that only the probabilities along the path from $[n]$ to $A$ remain. This yields the following result.

**Proposition 5.** *Let* $\mathbf{t} \in \mathbb{B}_{[n]}$ *and* $A \in \mathbf{t}$. *Denote by*

$$[n] = A_1 \supset \cdots \supset A_h = A$$

**Table 1.** Closed form expressions of the parameters for $\beta = -3/2, -1, 0, \infty$

| $\beta$ | $-3/2$ | $-1$ | $0$ | $\infty$ |
|---|---|---|---|---|
| $w(n)$ | $\dfrac{(2n-2)!}{2^{2n-2}(n-1)!}$ | $(n-1)!$ | $n!$ | $1$ |
| $Z(n)$ | $\dfrac{(2n-2)!}{2^{2n-3}(n-1)!}$ | $(n-1)!\sum_{j=1}^{n-1}\dfrac{1}{j}$ | $\dfrac{1}{2}(n-1)n!$ | $2^{n-1}-1$ |
| $\psi(n)$ | $\dfrac{1}{2}$ | $1\Big/\sum_{j=1}^{n-1}\dfrac{1}{j}$ | $\dfrac{2}{n-1}$ | $\dfrac{1}{2^{n-1}-1}.$ |

*the path from $[n]$ to $A$. We refer to $h \geq 1$ as the* height *of $A$ in* **t**. *The probability that $n+1$ attaches below $A$ is then*

$$\left(\prod_{j=1}^{h-1} \frac{p(\#A_{j+1}+1, \#(A_j \setminus A_{j+1}))}{p(\#A_{j+1}, \#(A_j \setminus A_{j+1}))}\right) p(\#A_h, 1).$$

For the uniform model (Gibbs fragmentation with $\beta = -3/2$), this product is telescoping, or we calculate directly from (8)

$$\left(\prod_{j=1}^{h-1} \frac{p(\#A_{j+1}+1, \#(A_j \setminus A_{j+1}))}{p(\#A_{j+1}, \#(A_j \setminus A_{j+1}))}\right) p(\#A_h, 1) = \frac{1}{2n-1},$$

giving a simple sequential construction (see, e.g., [15], Exercise 7.4.11).

It was shown in [9] that consistent Markovian fragmentation models can be assigned consistent independent exponential edge lengths, where the edge below vertex $A$ is given parameter $\lambda_{\#A}$, for a family $(\lambda_m)_{m \geq 1}$ of rates, where $\lambda_1 = 0$, $\lambda_2$ is arbitrary and $\lambda_m$, $m \geq 3$, is determined by $\lambda_2$ and the splitting rule $p$, in that consistency requires

$$\lambda_{n+1}(1 - p(n,1)) = \lambda_n \qquad \text{for all } n \geq 2. \tag{14}$$

The interpretation is that the partition of $[n+1]$ in $T_{n+1}$ (arriving at rate $\lambda_{n+1}$) splits $[n]$ only with probability $1 - p(n,1)$ and this thinning must reduce the rate for the partition of $[n]$ in $T_n$ to $\lambda_n$. This rate $\lambda_n$ also applies in $T_{n+1}$ after a first split $\{[n], \{n+1\}\}$.

Using consistency, equation (14) also implies

$$\lambda_n p(i,j) = \lambda_{n+1}(p(i,j+1) + p(i+1,j)) \qquad \text{for all } i, j \geq 1 \text{ with } i+j = n.$$

For the Gibbs fragmentation models, we obtain, using (14), (7), (12) and (13),

$$\lambda_n = \lambda_2 \prod_{j=2}^{n-1} \frac{1}{1 - p(j,1)} = \lambda_2 \prod_{j=2}^{n-1} \frac{Z(j+1)}{Z(j+1) - w(j)} = \lambda_2 Z(n) \prod_{j=2}^{n-1} \frac{1}{W_1 + W_{j-1}}$$

$$= \lambda_2 Z(n) \prod_{j=2}^{n-1} \frac{w(j-1)}{w(2)w(j-1) + w(j)} = \lambda_2 Z(n) \frac{\Gamma(4+2\beta)}{\Gamma(n+2+2\beta)},$$

where we require $\beta < \infty$ for the last step. Table 2 contains the rate sequences for $\beta = -3/2, -1, 0, \infty$ in the case $\lambda_2 = 1$.

Not only is $(\lambda_n)_{n \geq 3}$ determined by $p$, but a converse of this also holds.

**Proposition 6.** *Let $(\lambda_n)_{n \geq 2}$ be a consistent rate sequence associated with a consistent Markovian binary fragmentation model with splitting rule $p$, meaning that (14) holds. Then, $p$ is uniquely determined by $(\lambda_n)_{n \geq 2}$.*

**Proof.** It is evident from (14) that $p(n,1)$ is determined for all $n \geq 2$, and $p(1,1) = 1$. Now, (4) for $i = 1$ determines $p(i+1, j)$ for all $j \geq 2$, and an induction in $i$ completes the proof. □

A more subtle question is to ask what sequences $(\lambda_n)_{n\geq 2}$ arise as consistent rate sequences. The above argument can be made more explicit to yield

$$p(k, n-k) = \frac{1}{\lambda_n} \sum_{j=0}^{k} (-1)^{k-j+1} \binom{k}{j} \lambda_{n-j}, \qquad 1 \leq k \leq n/2,$$

which means that $(\lambda_n)_{n\geq 2}$ must have a discrete complete monotonicity, in that $k$th differences of $(\lambda_n)_{n\geq 2}$ must be of alternating signs, $k \geq 1$. This condition is not sufficient, however, as simple examples for $n = 3$ show ($\lambda_n = (n-1)^\alpha$ is completely monotone for $\alpha \in (0,1)$, but exchangeability implies that $1/3 = p(1,2) = (\lambda_3 - \lambda_2)/\lambda_3$ and so $\lambda_3 = 3/2$, whereas $(3-1)^\alpha \in (1,2)$ – even in the multifurcating case, cf. Section 5, we always have $\lambda_3 \leq 3/2$).

**Proposition 7.** *A sequence $(\lambda_n)_{n\geq 2}$ arises as rate sequence of a consistent Markovian binary fragmentation model if and only if*

$$\lambda_n = nc + \int_{(0,1)} (1 - x^n - (1-x)^n)\nu(\mathrm{d}x)$$

*for some $c \geq 0$ and $\nu$ a symmetric measure on $(0,1)$ with $\int_{(0,1)} x(1-x)\nu(\mathrm{d}x) < \infty$. The characteristics of the splitting rules associated with $(\lambda_n)_{n\geq 2}$ are $(c, \nu)$.*

**Proof.** This is a consequence of the integral representation (5) and [9], Proposition 3. Specifically, the association with Bertoin's theory of homogeneous fragmentations yields that each of $1, \ldots, n$ suffer erosion (being turned into a singleton) at rate $c$; the measure $\nu(\mathrm{d}x)$ gives the rate of fragmentations into two parts, to which $1, \ldots, n$ are allocated independently with probabilities $(x, 1-x)$, hence splitting $[n]$ with probability $1 - x^n - (1-x)^n$. □

The complete monotonicity is related to the study of the block containing 1, a *tagged fragment*; see [4, 10]. Since $\lambda_n$ is the rate at which one or more of $\{2, \ldots, n\}$ leave the block containing 1, the rate is composed of three components – a rate $c$ for the erosion of 1, a rate $(n-1)c$ for the erosion of $2,\ldots,n$ and a rate $\Lambda(\mathrm{d}z)$ of fragmentations into two parts, to which $2,\ldots,n$ are allocated independently with probabilities $(\mathrm{e}^{-z}, 1-\mathrm{e}^{-z})$, with 1 in the former part, hence splitting $[n]$ with probability $1 - \mathrm{e}^{-(n-1)z}$. Therefore

$$\lambda_n = c + (n-1)c + \int_{(0,\infty)} (1-\mathrm{e}^{-(n-1)z})\Lambda(\mathrm{d}z) = cn + \int_{(0,1)} \frac{1-\xi^{n-1}}{1-\xi}\mu(\mathrm{d}\xi) = \Phi(n-1)$$

for a Bernstein function $\Phi$, a finite measure $\mu$ on $(0,1)$ or a Lévy measure $\Lambda$ on $(0,\infty)$ with $\int_{(0,\infty)} (1 \wedge x)\Lambda(\mathrm{d}x) < \infty$; (see [4, 8, 10]), that is, $\lambda_n$ can be extended to a completely monotone function of a real parameter.

**Table 2.** Explicit rate sequences for $\beta = -3/2, -1, 0, \infty$

| $\beta$ | $-3/2$ | $-1$ | $0$ | $\infty$ |
|---|---|---|---|---|
| $\lambda_n$ | $\dfrac{n-1}{2^{2n-3}} \dbinom{2n-2}{n-1}$ | $\displaystyle\sum_{j=1}^{n-1} \frac{1}{j}$ | $\dfrac{3n-3}{n+1}$ | $2(1 - 2^{-(n-1)})$. |

## 5. Multifurcating Gibbs fragmentations and Poisson–Dirichlet models

As a generalization of the binary framework of the previous sections, we consider in this section consistent Markovian fragmentation models with splitting rule $p$ as in (2) of the Gibbs form

$$p(n_1,\ldots,n_k) = \frac{a(k)}{c(n)} \prod_{i=1}^{k} w(n_i) \tag{15}$$

for some $w(j) \geq 0$, $j \geq 1$, $a(k) \geq 0$, $k \geq 2$, and normalization constants $c(n) > 0$, $n \geq 2$. Note that we must have $w(1) > 0$ and $a(2) > 0$ to get positive probabilities for $n = 2$. To remove overparameterization, we will assume $w(1) = 1$ and $a(2) = 1$. Also, if we multiply $w(j)$ by $b^{j-1}$ and $a(k)$ by $b^k$ (and $c(n)$ by $b^n$), the model remains unchanged. We will use this observation to get a nice parameterization in the consistent case (Theorem 8 below).

In [9], we showed that consistency of the model is equivalent to the set of equations

$$p(n_1,\ldots,n_k) = p(n_1+1, n_2,\ldots,n_k) + \cdots + p(n_1,\ldots,n_k+1) + p(n_1,\ldots,n_k,1) \\ + p(n_1+\cdots+n_k, 1)p(n_1,\ldots,n_k) \tag{16}$$

for all $n_1,\ldots,n_k \geq 1$, $k \geq 2$. We also established an integral representation extending (5) to the multifurcating case. The special case relevant for us is in terms of a measure $\overline{\nu}$ on $\mathcal{S}^\downarrow = \{s = (s_i)_{i\geq 1} : s_1 \geq s_2 \geq \cdots \geq 0, s_1 + s_2 + \cdots = 1\}$ satisfying $\int_{\mathcal{S}^\downarrow} (1-s_1)\overline{\nu}(\mathrm{d}s) < \infty$:

$$p(n_1,\ldots,n_k) = \frac{1}{Z(n_1+\cdots+n_k)} \int_{\mathcal{S}^\downarrow} \sum_{i_1,\ldots,i_k \text{ distinct}} \prod_{j=1}^{k} s_{i_j}^{n_j} \overline{\nu}(\mathrm{d}s). \tag{17}$$

The general case has a further parameter $c \geq 0$, as in (5), and also allows $\overline{\nu}$ to charge $(s_i)_{i\geq 1}$ with $s_1 + s_2 + \cdots < 1$; see [9]. We will only meet the extreme case $p(1,\ldots,1) = 1$, which corresponds to $\overline{\nu} = \delta_{(0,0,\ldots)}$.

We set

$$\frac{a(k+1)}{a(k)} = A_k, \qquad \frac{c(n+1)}{c(n)} = C_n, \qquad \frac{w(n+1)}{w(n)} = W_n$$

and, in analogy to Proposition 5, we find that, given $T_n = \mathbf{t} \in \mathbb{T}_{[n]}$, for each vertex $B \in \mathbf{t}$, the probability that $n+1$ attaches below $B$ is

$$\left(\prod_{j=1}^{h-1} \frac{W_{n_{j+1}}}{C_{n_j}}\right) \frac{a(2)w(n_h)w(1)}{c(n_h+1)},$$

where $[n] \supset S_1 \supset \cdots \supset S_h = B$ is the path from $[n]$ to $B$, $n_j = \#S_j$ and $k_j$ denotes the number of children of $S_j$, $j = 1, \ldots, h$.

However, $n+1$ can also attach as a singleton block to an existing partition $\{B_1, \ldots, B_k\}$ of $B \in T_n$. In this case, we say that $n+1$ attaches to the vertex $B$. For each non-leaf vertex $B \in \mathbf{t}$, the probability that $n+1$ attaches to the vertex $B$ is

$$\left(\prod_{j=1}^{h-1} \frac{W_{n_{j+1}}}{C_{n_j}}\right) \frac{A_{k_h} w(1)}{C_{n_h}}.$$

In this framework, we have the following generalization of Theorem 2 to the multifurcating case.

**Theorem 8.** *If $p$ is of the Gibbs form (15) and consistent, then $p$ is associated with the two-parameter Ewens–Pitman family given by*

$$w(n) = \frac{\Gamma(n-\alpha)}{\Gamma(1-\alpha)}, \qquad n \geq 1, \quad \text{and} \quad a(k) = \alpha^{k-2} \frac{\Gamma(k+\theta/\alpha)}{\Gamma(2+\theta/\alpha)}, \qquad k \geq 2$$

*(or limiting quantities $\alpha \downarrow 0$), $c(n)$, $n \geq 1$, being normalization constants, for a parameter range extended as follows:*

- *either $0 \leq \alpha < 1$ and $\theta > -2\alpha$ (multifurcating cases with arbitrarily high block numbers),*
- *or $\alpha < 0$ and $\theta = -m\alpha$ for some integer $m \geq 3$ (multifurcating with at most $m$ blocks),*
- *or $\alpha < 1$ and $\theta = -2\alpha$ (binary case),*
- *or $\alpha = -\infty$ and $\theta = m$ for some integer $m \geq 2$, that is, $a(2) = 1$, $a(k) = (m-2)\cdots(m-k+1)$, $k \geq 3$, and $w(j) \equiv 1$ (recursive coupon collector, where a split of $[n]$ is obtained by letting each element of $[n]$ pick one of $m$ coupons at random, just conditioned so that at least two different coupons are picked),*
- *or $\alpha = 1$, that is, $w(1) = 1$, $w(j) = 0$, $j \geq 2$ (deterministic split into singleton blocks).*

*In terms of the integral representation (17), the measure $\overline{\nu}$ on $\mathcal{S}^{\downarrow}$ is, respectively, size-ordered Poisson–Dirichlet$(\alpha, \theta)$, Dirichlet$(-\alpha, \ldots, -\alpha)$, Beta$(-\alpha, -\alpha)$, $\delta_{(1/m, \ldots, 1/m)}$ and $\delta_{(0,0,\ldots)}$.*

**Proof.** For the Gibbs fragmentation model with $w(1) = a(2) = 1$ and $w(j) > 0$ for all $j \geq 2$ with notation as introduced, consistency (16) is easily seen to be equivalent to

$$C_n = W_{n_1} + \cdots + W_{n_k} + A_k + \frac{w(n)}{c(n)} \qquad \text{for all } n_1 + \cdots + n_k = n, \tag{18}$$

where $k \leq m$ if $m = \inf\{i \geq 1 : a(i+1) = 0\} < \infty$.

As in the proof of Theorem 2, we deduce from this (the special case $k = 2$) that either $W_j = a > 0$ (excluded for the time being as $b = 0$) or

$$W_j = a + bj \quad \Rightarrow \quad w(j) = W_1 \ldots W_{j-1} = b^{j-1} \frac{\Gamma(j-\alpha)}{\Gamma(1-\alpha)} \qquad \text{for all } j \geq 1,$$

for some $b > 0$, $a > -b$ and $\alpha := -a/b < 1$. As noted above, we can reparameterize so that we get $b = 1$ without loss of generality. In particular, $W_j = j - \alpha$, $j \geq 1$, and so (18) reduces to

$$C_n = n - k\alpha + A_k + \frac{w(n)}{c(n)} \qquad \text{for all } 2 \leq k \leq m \wedge n.$$

Similarly, we deduce that $\theta := A_k - k\alpha$ does not depend on $k$ and so $a(k) = \theta^{k-2}$ if $\alpha = 0$, and otherwise,

$$A_k = \theta + k\alpha \quad \Rightarrow \quad a(k) = A_2 \ldots A_{k-1} = \alpha^{k-2} \frac{\Gamma(k+\theta/\alpha)}{\Gamma(2+\theta/\alpha)} \qquad \text{for all } 2 \leq k \leq m+1.$$

Note that this algebraic derivation leads to probabilities in (15) only in the following cases.

- If $0 \leq \alpha < 1$, then $a(3) = A_2 = \theta + 2\alpha > 0$ if and only if $\theta > -2\alpha$, and then also $A_k = \theta + k\alpha > 0$ and $a(k) > 0$ for all $k \geq 3$.
- If $\alpha < 0$, then $a(3) = A_2 = \theta + 2\alpha > 0$ if and only if $\theta > -2\alpha$ also, but then $A_k = \theta + k\alpha$ is strictly decreasing in $k$ and $A_k < 0$ eventually, which impedes $m = \infty$. If we have $m < \infty$, we achieve $a(m+1) = 0$ if and only if $\theta = -m\alpha$. The iteration only takes us to $a(m+1) = 0$ and we specify $a(k) = 0$ for $k > m$ also. We cannot specify $a(k)$, $k > m+1$, differently, since every consistent Gibbs fragmentation with $a(k) > 0$ for $k > m+1$ has the property that $T_{[k]} = \{[k], \{1\}, \ldots, \{k\}\}$ has only one branch point $[k]$ of multiplicity $k$ with positive probability, but then the restricted tree $T_{[m+1],[k]} = \{[m+1], \{1\}, \ldots, \{m+1\}\}$ with positive probability, which contradicts $a(m+1) = 0$.
- If $a(3) = 0$, that is, $m = 2$, the argument of the preceding bullet point shows that we are in the binary case $a(k) = 0$ for all $k \geq 3$ and we can conclude by Theorem 2.
- The case $b = 0$ is the limiting case $\alpha = -\infty$ with $w(j) \equiv 1$. We take up the argument to see that $A_k = \theta - k$ and so $m < \infty$ and $\theta = m$, where we then get $a(2) = 1$ and $a(k) = (m-2) \cdots (m-k+1)$, $3 \leq k \leq m+1$.

Finally, if $w(m) = 0$ for some $m \geq 2$, then consistency imposes $w(j) = 0$ for all $j \geq m$, and it follows from the integral representation (17) that in fact $w(j) = 0$ for all $j \geq 2$. The identification of $\overline{\nu}$ on the standard parameter range can be read from [15], Section 3.2. For the extension $-\alpha \geq \theta \geq -2\alpha$, we refer to [10]. $\square$

Kerov [11] showed that the only exchangeable partitions of $\mathbb{N}$ of Gibbs type are of the two-parameter family $\text{PD}(\alpha, \theta)$ with usual range for parameters $\theta > -\alpha$, etc.; see also [7, 14]. Theorem 8 is a generalization to splitting rules that allows an extended parameter range for the same reason as in the binary case: the trivial partition of one single block is excluded from $p$ and when associating consistent exponential edge lengths with parameters $\lambda_m$, $m \geq 1$, the first split of $[m+1]$ happens at a higher and higher rate and we may have $\lambda_m \to \infty$. In fact,

$$\kappa(\{\pi \in \mathcal{P}_\mathbb{N} : \pi|_{[n]} = \{B_1, \ldots, B_k\}\}) = \lambda_n p(\#B_1, \ldots, \#B_k)$$

uniquely defines a $\sigma$-finite measure on $\mathcal{P}_\mathbb{N} \setminus \{\mathbb{N}\}$, the set of non-trivial partitions of $\mathbb{N}$, associated with a homogeneous fragmentation process. This is closely related to (17) via Kingman's paintbox representation $\kappa = \int_{\mathcal{S}^{\downarrow}} \kappa_s \overline{\nu}(\mathrm{d}s)$. The extended range was first observed by Miermont [13] in the special case $\theta = -1$ (related to the stable trees of Duquesne and Le Gall [5]).

We refer to [10] for a study of spinal partitions of Markovian fragmentation models. There are notions of fine and coarse spinal partitions. First, remove from $T_n$ the spine of 1, that is, the path from $[n]$ to $\{1\}$. The resulting collection is a disjoint union of fragmentations of sets $B_j$, say, that form a partition of $\{2, \ldots, n\}$, which is called the *fine spinal partition*. Second, merge blocks (in the multifurcating case) that were children of the same spinal vertex; the resulting partition is called the *coarse spinal partition*. It is shown that for the splitting rules from the two-parameter family with parameters $\alpha$ and $\theta$ (the Gibbs fragmentations), the fine partition is obtained from the coarse partition by applying independently for each block of the coarse partition an exchangeable partition from the two-parameter family of random partitions, with parameters $\alpha$ and $\alpha + \theta$.

## Acknowledgements

This research was supported in part by EPSRC Grant GR/T26368/01 and NSF Grants DMS-04-05779 and DMS-03-05009. M. Winkel was also supported by the Institute of Actuaries and the insurance group Aon Limited.

## References

[1] Aldous, D. (1991). The continuum random tree. I. *Ann. Probab.* **19** 1–28. MR1085326
[2] Aldous, D. (1996). Probability distributions on cladograms. In *Random Discrete Structures (Minneapolis, MN, 1993). IMA Vol. Math. Appl.* **76** 1–18. New York: Springer. MR1395604

[3] Berestycki, N. and Pitman, J. (2007). Gibbs distributions for random partitions generated by a fragmentation process. *J. Stat. Phys.* **127** 381–418. MR2314353

[4] Bertoin, J. (2001). Homogeneous fragmentation processes. *Probab. Theory Related Fields* **121** 301–318. MR1867425

[5] Duquesne, T. and Le Gall, J.-F. (2002). Random trees, Lévy processes and spatial branching processes. *Astérisque* **281** vi+147. MR1954148

[6] Ford, D.J. (2005). Probabilities on cladograms: Introduction to the alpha model. Preprint. arXiv:math.PR/0511246.

[7] Gnedin, A. and Pitman, J. (2005). Exchangeable Gibbs partitions and Stirling triangles. *Zap. Nauchn. Sem. S.-Peterburg. Otdel. Mat. Inst. Steklov. (POMI)* **325** (*Teor. Predst. Din. Sist. Komb. i Algoritm. Metody* **12**) 83–102, 244–245. MR2160320

[8] Gnedin, A. and Pitman, J. (2006). Moments of convex distribution functions and completely alternating sequences. Preprint. arXiv:math.PR/0602091.

[9] Haas, B., Miermont, G., Pitman, J. and Winkel, M. (2006). Continuum tree asymptotics of discrete fragmentations and applications to phylogenetic models. Preprint. arXiv:math.PR/0604350. *Ann. Probab.* To appear.

[10] Haas, B., Pitman, J. and Winkel, M. (2007). Spinal partitions and invariance under re-rooting of continuum random trees. Preprint. arXiv:0705.3602. *Ann. Probab.* To appear.

[11] Kerov, S. (2005). Coherent random allocations, and the Ewens–Pitman formula. *Zap. Nauchn. Sem. S.-Peterburg. Otdel. Mat. Inst. Steklov. (POMI)* **325** (*Teor. Predst. Din. Sist. Komb. i Algoritm. Metody* **12**) 127–145, 246. MR2160323

[12] McCullagh, P., Pitman, J. and Winkel, M. (2007). Gibbs fragmentation trees. Preprint. arXiv:0704.0945.

[13] Miermont, G. (2003). Self-similar fragmentations derived from the stable tree. I. Splitting at heights. *Probab. Theory Related Fields* **127** 423–454. MR2018924

[14] Pitman, J. (2003). Poisson–Kingman partitions. In *Statistics and Science: A Festschrift for Terry Speed. IMS Lecture Notes Monogr. Ser.* **40** 1–34. Beachwood, OH: Inst. Math. Statist. MR2004330

[15] Pitman, J. (2006). *Combinatorial Stochastic Processes. Lecture Notes in Math.* **1875**. *Lectures from the 32nd Summer School on Probability Theory held in Saint-Flour, July 7–24, 2002.* Berlin: Springer. MR2245368

[16] Schroeder, E. (1870). Vier combinatorische Probleme. *Z. f. Math. Phys.* **15** 361–376.

[17] Semple, C. and Steel, M. (2003). *Phylogenetics. Oxford Lecture Series in Mathematics and Its Applications* **24**. Oxford Univ. Press. MR2060009

[18] Stanley, R.P. (1999). *Enumerative Combinatorics.* **2**. *Cambridge Studies in Advanced Mathematics* **62**. Cambridge Univ. Press. MR1676282